\documentclass{article}

\usepackage{latexsym,amssymb}

\title{Bounding the order of finite $p^{\prime}$- subgroups of ${\rm GL}(n,\mathbb{C})$}

\author{Geoffrey R. Robinson}

\begin{document}

\maketitle

\begin{abstract}

In this note, which originated in response to a question from P. Etingof (related to a potential application to the forthcoming paper [3]), we prove the following:

\medskip
\noindent {\bf Theorem  A:} \emph{ There is a fixed constant $C$ such that for any positive integer $n$ and prime $p$, every finite subgroup $G$ of order coprime to $p$ of  
${\rm GL}(n,\mathbb{C})$ has an Abelian normal subgroup  $A$ with $$[G:A] \leq (Cp)^{n-1}.$$}

\medskip
This theorem may be regarded as a refinement of  the classical theorem of Jordan to finite complex linear groups $G$ of order not divisible by a specified prime 
$p$  (though when $p > \frac{n}{eC}$ it offers no real improvement on the generically worst case bound $(n+1)!$ for Jordan's theorem).
For any given prime $p$, though, this bound will be much less than $(n+1)!$ for sufficiently large $n$.

\medskip
Consideration of the symmetric group $Sym_{p-1}$ shows that whenever $0 < c < \frac{1}{e}$, for any sufficiently large prime $p$, there is a $p^{\prime}$-subgroup $G$ of ${\rm GL}(p-2,\mathbb{C})$ with no non-identity Abelian normal subgroup, yet with $|G| >(cp)^{p-3}.$ More generally, for any such $c$, and all sufficiently large $p$,  we can find (for infinitely many values of $n$) a $p^{\prime}$-subgroup $G$ of ${\rm GL}(n,\mathbb{C})$ with no non-identity Abelian normal subgroup, yet with 
$|G| >(cp)^{n-1}.$

\medskip
Once Theorem A is proved, a slightly more general version (with a different constant) follows easily. We recall that when $p$ is a prime number, the finite group $G$
is said to be $p$-\emph{closed} if $G$ has a unique (necessarily normal) Sylow $p$-subgroup.

\medskip
\noindent {\bf Theorem  B:} \emph{ There is a fixed constant $C^{\ast}$ such that for any positive integer $n$ and prime $p$, every finite $p$-closed  subgroup $G$ of  
${\rm GL}(n,\mathbb{C})$ has an Abelian normal subgroup  $A^{\ast}$ with $$[G:A^{\ast}] \leq (C^{\ast}p)^{n-1}.$$ Furthermore, if the Sylow $p$-subgroup of the $p$-closed group $G$ is Abelian, then we may choose $A^{\ast}$ so that $[G:A^{\ast}]$ is coprime to $p$.}

\end{abstract}

\section{Preliminary Discussion}

\medskip
Let $G$ be a finite subgroup of ${\rm GL}(n,\mathbb{C})$. Then there is an Abelian normal subgroup $A$ of $G$ with $[G:A]$ bounded in terms of $n$ alone, a classical theorem of 
C. Jordan. Allowing the presence of an Abelian normal subgroup is unavoidable. Explicit upper bounds for $[G:A]$ have been given by various authors,  but since the symmetric group $Sym_{n+1}$ has a faithful complex representation of degree $n$ for $n \geq 2$, this sets a lower bound of $(n+1)!$ for an optimal explicit bound. Prior to his unfortunate disappearance, B. Weisfeiler outlined a strategy to use the Classification of Finite Simple Groups to obtain substantial improvements on the previously known explicit bounds. In [2], M.J. Collins later proved that we do have  $[G:A] \leq (n+1)!$ when  $n \geq 71$, and explicitly determined the optimal upper bounds for $n < 71$.

\medskip
Our purpose here is to show that if we restrict the prime divisors of $|G|$, then we may often improve these bounds for the general finite complex linear group $G$. When $p = n+2$ is  a prime greater than $71$, we note that ${\rm Sym}_{n+1}$ is a subgroup of order prime to $p$ of ${\rm GL}(n,\mathbb{C})$, and the worst case bound is attained, so such improvement is not always possible. As noted in the abstract, we were led to consider the case of general finite complex linear groups of order not divisible by an arbitrary chosen prime $p$ because of a question from P. Etingof, in connection with $[3]$.

\medskip
In case $G$ is solvable, it is possible to improve the general bound for Jordan's theorem significantly. By a theorem of L. Dornhoff, [5], there is a constant $c$ such that $$[G:A] \leq c^{n-1}$$ for some Abelian normal subgroup $A$ of $G$. In particular, we have $$[G:A] < 2^{\frac{4n}{3}}3^{\frac{10n}{9}}$$ for a suitable Abelian normal subgroup $A$ of $G$. A key ingredient in the proof of this result is the fact that a solvable subgroup of the symmetric group $Sym_{n}$ has order at most $24^{\frac{n-1}{3}}$, a bound which may be attained  whenever $n = 4^{k}$ for a positive integer $k$. In [6], R. Gow showed that a  smaller constant $c^{\prime}$ suffices for finite complex linear groups $G$ of odd order. \medskip
In [11], the present author improved the bounds of Dornhoff and of Gow for finite complex linear groups of odd  order,  obtaining bounds which improved as the smallest prime dividing the group order increased in that case.  

\medskip
If $p$ is a prime, and $G$ is a $p^{\prime}$-subgroup of $Sym_{n}$, that is a subgroup of order prime to $p$ of $Sym_{n}$ ($n \geq 2$), then P\'alfy and Pyber proved in [9]  that  
$|G| < p^{n-1}.$  On the other hand, the $d$-fold iterated wreath product $${\rm Sym}_{p-1} \wr Sym_{p-1} \wr \ldots \wr Sym_{p-1}$$ illustrates that whenever  
$c$ is a real number strictly less than $\frac{1}{e}$, then for each sufficiently large prime $p$, the symmetric group $Sym_{n}$ has a subgroup $G$ of order greater than $(cp)^{n-1}$ whenever $n =(p-1)^{d}$ for a positive integer $d$. Note also that such a group $G$ has no non-identity Abelian normal subgroup, and is isomorphic to a subgroup of 
${\rm GL}(n,\mathbb{C}).$
																																																											\medskip
We will prove:

\medskip
\noindent {\bf Theorem A:} \emph{ There is a constant $C$ (independent of both the chosen integer $n$ and the chosen prime $p$) such that whenever $p$ is a prime and $G$ is a finite $p^{\prime}$-subgroup of 
${\rm GL}(n,\mathbb{C})$, then $G$ has an Abelian normal subgroup $A$  with $$[G:A] \leq (Cp)^{n-1}.$$}

\medskip
We make no attempt to find the optimal choice of constant $C$, however, and we work with very crude estimates in order to facilitate a smooth proof. The fact that ${\rm SL}(2,5)$ is isomorphic to a subgroup $G$ of ${\rm GL}(2,\mathbb{C})$  illustrates ( with $n =2$ and $p=7)$ that we certainly need $$C \geq  \frac{60}{7} > 8.57.$$

\medskip
Our proof requires the Classification of the Finite (non-Abelian) Simple Groups: the main properties we will use are the following, ( i)-iii) may be checked using [4] : for iv), it is well known that alternating groups may be generated by two elements, and it is a theorem of R. Steinberg that simple groups of Lie type may be generated by two elements).

\medskip
\noindent i) There are only finitely many isomorphism types of non-Abelian finite simple groups which are neither alternating nor of Lie type, (and  each such exceptional (sporadic) simple group has outer automorphism group of order at most two, though this latter fact is not really essential to our proof-the fact that $|{\rm Out}(S))|$ is bounded for any sporadic simple group would suffice).

\medskip
\noindent ii) The outer automorphism group of the alternating group $Alt_{d}$ has order two for $d > 6$ and for $d = 5$. The outer automorphism group of $Alt_{6}$ is 
a Klein $4$-group.

\medskip
\noindent iii) Each finite simple group $S$ of Lie type of characteristic $\ell$ has a nilpotent subgroup $U$ of order greater than $|S|^{\frac{1}{3}}.$ This is a slight extension of an observation of E. Artin.

\medskip
\noindent iv) If $S$ is a non-Abelian finite simple group, then $S$ may be generated by two elements, and consequently (on consideration of the image of a pair of generators under an automorphism), we have $|{\rm Aut}(S)| < |S|^{2}.$ Notice that this implies that $|{\rm Out}(S)| < |S|.$ If $S$ is sporadic or alternating of degree other than $6$, we have 
$|{\rm Out}(S)| \leq 2.$

\medskip
Given these facts, it is possible to write down an explicit choice of constant $C$, though our proof in its current form requires a choice of $C$ which is far from optimal.

\medskip
We outline the structure of the proof now:  We first show that it is sufficient to find a constant $D$ (independent of both $n$ and $p$) such that $$[G:F(G)] \leq (Dp)^{n-1},$$ where
$F(G)$ denotes the Fitting subgroup (the unique largest nilpotent normal subgroup) of $G$. This is accomplished in the first section,  where we use a careful analysis of classical results of Blichfeldt and Frobenius (see, for example [1], [5],[7],[10] and [12]), to demonstrate that if $H$ is a 
finite subgroup of ${\rm GL}(m,\mathbb{C})$,  then $H$ has an Abelian normal subgroup $A$ such that any nilpotent subgroup $N$ of $H$ has
$[N: N \cap A] \leq (24)^{m-1}$ (so that, in particular, we have $A \leq F(H)$ and $[F(H):A] \leq (24)^{m-1}).$

\medskip
In the second section, we use induction on $n$ to reduce to the case that $G$ is irreducible and primitive.

\medskip
In the third section, we choose a non-central  normal subgroup $M$ of $G$ which is minimal subject to strictly containing $Z(G)$.
At this point, there are two scenarios to investigate. We either have $M = F(M)$, or else $M = Z(G)E(M)$, where $E(M) = L_{1}L_{2} \ldots L_{t}$, with the $L_{i}$ being $G$-conjugate components (quasi-simple subnormal subgroups) of $M$ (possibly $t =1$).

\medskip
In each of these cases, we need to consider the structure and order of $${\rm Out}_{G}(M) \cong G/MC_{G}(M),$$ and the proof subdivides further, according to whether $C_{G}(M) \leq M$ or not.  If $C_{G}(M) \leq M$, then we see that $M = F^{\ast}(G),$ the generalized Fitting subgroup of $G$. If, in that case, $M = F(M)$, then we see (by minimality of $M$), that\\ $M = Z(G)U$ with $U$ extra-special (or $U = \Omega_{2}(O_{q}(M))$ when $q=2$), and that ${\rm Out}_{G}(M)$  is isomorphic to  a subgroup of the symplectic group ${\rm Out}(U).$  In that case, there is a suitable choice of constant $D$. 

\medskip
On the other hand, (still in the case $C_{G}(M) = Z(M)),$ if we have\\ $M = Z(G)E(M),$ then a detailed examination of the case supplies a suitable constant $D$ (which may be chosen large enough to also cover the case\\ $M = F(M)$).

\medskip
Having established the existence of a suitable constant $D$ to cover the case $C_{G}(M) = Z(M)$, we then prove by induction on the dimension $n$ that the same constant $D$
works in general. Notice that when $C_{G}(M) \neq Z(M)$ then both $M$ and $C_{G}(M)$ are isomorphic to complex linear groups in dimension strictly less than $n$,while
the structure of $$G/MC_{G}(M) \cong {\rm Out}_{G}(M)$$ is as previously analyzed. 

\medskip
Having outlined the structure of the proof, we now commence the proof proper.

\section{Nilpotent subgroups of finite complex linear groups}

\medskip
Let $H$ be a finite subgroup of ${\rm GL}(m,\mathbb{C})$.  It follows from arguments due to Frobenius that $H$ has an Abelian normal subgroup $A$ which is the subgroup generated
by those elements of $H$ whose eigenvalues all lie on arc of length at most $\frac{\pi}{3}$ on the unit circle $S^{1}.$ This is explained in some detail in [10]. An argument
originating with Frobenius shows that each Abelian subgroup $B$ of $H$ satisfies $[B: B \cap A] \leq 12^{m-1}.$ 

\medskip
Now let $N$ be a nilpotent subgroup of $H$. Then $N$ is conjugate to a monomial subgroup of ${\rm GL}(m,\mathbb{C})$, and that monomial group has an Abelian normal subgroup (consisting of its diagonal matrices) with factor group isomorphic to a nilpotent subgroup of the symmetric group ${\rm Sym}_{m}.$ It is well-known and easy to
check that each nilpotent subgroup of $Sym_{m}$ has order at most $2^{m-1}$  (a bound which can be attained whenever $m$ is a power of $2$). Hence $N$ has an Abelian normal subgroup $T$ with $[N:T] \leq 2^{m-1}$ , and we see that $$[N:T \cap A] \leq (24)^{m-1}.$$ In particular, $[N: N \cap A] \leq (24)^{m-1}.$

\medskip
Now, considering the case $N = F(H)$  we have $A \leq F(H)$ and $$[F(H) : A] \leq (24)^{m-1}.$$ In particular, if we have $[H:F(H)] \leq D^{m-1}$ for some constant $D$, then we have $$[H : A] \leq (24D)^{m-1}.$$

\medskip
\noindent{\bf Remark:} Hence, returning to our finite subgroup $G$ of ${\rm GL}(n,\mathbb{C})$, if we can prove that there is a constant $D$, independent of $p$, such that 
$[G:F(G)] \leq (Dp)^{n-1}$ whenever $G$ has order prime to $p$, then $G$ has an Abelian normal subgroup of index at most $(Cp)^{n-1}$, where $C = 24D.$ Thus 
we have now reduced to the task of finding a constant $D$, independent of ($n$ and) the prime $p$, such that $[G:F(G)] \leq (Dp)^{n-1}$ whenever $G$ is a finite subgroup of 
${\rm GL}(n,\mathbb{C})$ of  order prime to $p$.

\medskip
\section{ Reduction to the primitive case}

\medskip
We now aim to prove by induction that there is a constant $D > 1$, independent of both $n$ and $p$, such that $$[G:F(G)] \leq (Dp)^{n-1}$$ whenever $G$ is a finite subgroup of
${\rm GL}(n,\mathbb{C})$ whenever $G$ is a finite subgroup of ${\rm GL}(n,\mathbb{C})$ of order prime to $p$, the result being clear when $n = 1$.

\medskip
Suppose that we have found a constant $D$ which works (for all primes $p$) for irreducible primitive finite complex linear groups in all degrees, 
and works (for all primes $p$) for all finite complex linear groups of degree less than $n$.   We claim that the constant $D$ works for our group $G$.

\medskip
If $G$ is reducible, then $G$ is isomorphic to a subgroup of a $p^{\prime}$-group $$G_{1} \times G_{2},$$ where there are positive integers $s,t $ with $s +t = n$,
and each $G_{i}$ is a homomorphic image of $G$, where $$G_{1} \leq {\rm GL}(s,\mathbb{C})$$ and $$G_{2} \leq {\rm GL}(t,\mathbb{C}).$$

\medskip
Then we have $$F(G_{1} \times G_{2}) = F(G_{1}) \times F(G_{2}),$$ and  $$[G_{1} \times G_{2}: F(G_{1} \times G_{2}) ] \leq (Dp)^{s-1}(Dp)^{t-1}
\leq (Dp)^{n-2}.$$ Hence $$[G:F(G)] \leq [G : G \cap [F(G_{1}) \times F(G_{2})] ]\leq (Dp)^{n-2},$$ and $D$ works for $G$.

\medskip
Suppose then that $G$ is irreducible, but imprimitive. Suppose then that the given representation of $G$ is equivalent to a representation induced from a primitive representation of a proper subgroup $H$ of $G$, say of degree $m$ (a divisor of $n$) . Let $t =[G:H]$. Let $N = \bigcap_{g \in G}H^{g} \lhd G$. Then $G/N$ is isomorphic to a $p^{\prime}$-subgroup of the symmetric group $Sym_{t}.$ By [9], we have $$[G : N] \leq p^{t-1}.$$ By Clifford's Theorem, $N$ embeds in a direct product 
$N_{1} \times N_{2} \ldots \times N_{t}$, where each $N_{i}$ is isomorphic to a finite (not necessarily irreducible) subgroup of ${\rm GL}(m,\mathbb{C}).$

\medskip
By the inductive hypothesis, we have $$[N_{i}:F(N_{i})] \leq (Dp)^{m-1}$$ for each $i$, and we then obtain $$[N:F(N)] \leq (Dp)^{tm-t}.$$ Hence
$$[G: F(N)] \leq (Dp)^{tm-t}p^{t-1} \leq (Dp)^{n-1}$$ and certainly $$[G:F(G)] \leq (Dp)^{n-1}.$$

\medskip
Hence it suffices to find a constant $D >1$, independent of $n$ and $p$, such that $$[G:F(G)] \leq (Dp)^{n-1}$$ whenever $G$ is a primitive (irreducible) 
finite subgroup of ${\rm GL}(n,\mathbb{C}).$

\medskip
\section{The primitive case: A minimal non-central normal subgroup of $G$}

\medskip
We now suppose that $G$ is a non-Abelian finite primitive (irreducible) $p^{\prime}$- subgroup of ${\rm GL}(n, \mathbb{C})$, (where $n >1$). In this case, we note that all 
Abelian normal subgroups of $G$ are central, and, more generally, (using Clifford's Theorem), that any normal subgroup $N$ of $G$ is isomorphic to an irreducible
subgroup of ${\rm GL}(d,\mathbb{C})$ for some divisor $d$ of $N$, while the normal subgroup $C_{G}(N)$ is isomorphic to a (not necessarily irreducible) subgroup 
of ${\rm GL}(\frac{n}{d},\mathbb{C}).$ 

\medskip
We may choose a normal subgroup $M$ of $G$ with $Z(G) \leq M$ and with $M$  minimal subject to being non-central in $G$. Then $M$ is non-Abelian by the primitivity of $G$, so that 
$$F^{\ast}(M) \neq Z(M) ( = Z(G)).$$ Since $F(M) \lhd G$ and $E(M) \lhd G$, the minimal choice of $M$  forces either\\ $M = F(M)$ or $M = Z(G) E(M).$ It may happen that $M = G$, so we consider this case first.

\medskip
In that case, since we are done if $G$ is nilpotent, see see that $G = Z(G)L$ with $L$ quasi-simple.

\medskip
More generally, for later convenience, we deal first with the case that that $G$ is a finite (irreducible) primitive subgroup of ${\rm GL}(n,\mathbb{C})$ with $F^{\ast}(G) = Z(G)L$, where $L = L^{\prime}$ and  $L/Z(L) \cong S,$  a finite non-Abelian simple group. Then $[G:F^{\ast}(G)] \leq |{\rm Out}(S)| $ and $[F^{\ast}(G): Z(G)] = |S|.$ In this case, we have\\ $F(G) = Z(G)$ and  $[G:F(G)] \leq |{\rm Aut}(S)|.$ Also, for later covenience, we note that there is a bijection between nilpotent subgroups of $S$ and nilpotent subgroups of $L$ which contain $Z(L)$.

\medskip
Suppose first that $S$ is sporadic. Since there are only finitely many isomorphism types of sporadic simple groups, and we have $|{\rm Out}(S)| \leq 2$ whenever $S$ is sporadic,
we may choose  a constant $D_{1}$ so that $D_{1} > |{\rm Aut}(S)|$, whichever sporadic group $S$ is, and then certainly $$[G:F(G)] \leq (D_{1}p)^{n-1}.$$ For example, we may take $D_{1}$ to be order of the Fischer-Griess Monster sporadic group $\mathcal{M}$. In fact, since any faithful complex representation of any perfect central extension of a sporadic simple group $S$ has dimension at least $6$, and each such $S$ has order at most  $|\mathcal{M}| < 10^{54},$ we see that $D_{1} = 10^{11}$ works.

\medskip
Suppose next that $S$ is simple of Lie type. Then by the remark in the preliminary discussion about Artin's observation, and the discussion of nilpotent subgroups of finite complex linear groups in the first section, we know that $S$ has a nilpotent subgroup $U$ with $|U|^{3} > |S|,$ and that $|U| \leq  (24)^{n-1}$ (since all Abelian normal subgroups of $G$ are central). Hence  we have $|S| \leq (24)^{3(n-1)}$. Also, we certainly have $|{\rm Out}(S)|<|S|$, so that $$[G:F(G)] \leq 24^{6(n-1)}.$$ In this case, we set 
$D_{2} = (24)^{6} < 10^{9}.$

\medskip
Suppose next that $S$ is an alternating group, say $S \cong Alt_{k}$ for some $k \geq 5$. Now we need to take account of the fact that $G$ has order prime to $p$.
Notice then that we certainly have $k < p$ since $G$ has order prime to $p$. Hence $p \geq 7.$ Then $$[G:Z(G)] \leq 2k! < 2k^{k}.$$ 

\medskip
It is well-known that for $k \geq 7$, the alternating group  $Alt_{k}$ has no faithful complex irreducible representation of degree less than $k-1$ 
By a Theorem of A.Wagner [13], for all but a small finite number of $k$, no perfect central extension of $Alt_{k}$ has a faithful complex representation  of degree
less than $k-1$. To err on the conservative side, whenever $k \geq 23$, no perfect central extension of $Alt_{k}$ has a faithful complex representation of degree less than $k-1$. Hence there is a constant $D_{3}$ (independent of $p$), such that $D_{3}^{n-1} \geq  4|Alt_{k}| \geq [G:Z(G)]$ in any case that $n < k-1$. 

\medskip
We now consider the case that $n \geq k-1$ and $n-1 \geq k-2$. 

\medskip
We now  show that, in that case, we may choose a constant $D_{4}$ (independent of $p$) so that $(D_{4}p)^{k-2} \geq 2k!$.

\medskip
But we certainly have $4^{k-2}  \geq (k-1)^{2}$ (as $k \geq 5$) and $p^{k-2} \geq (k+1)^{k-2}$. Thus  we always have 
$$(4p)^{k-2} \geq (k-1)^{2}(k+1)^{k-2}  \geq (k^{2}-1) (k-1)(k+1)^{k-3}> k!(k-1) > 4k! ,$$ 
so that $D_{4} = 4$ will work.

\medskip
In conclusion, in all cases where $F^{\ast}(G) = Z(G)L$ with $L$ quasi-simple, we may find a constant $D= max\{D_{1},D_{2},D_{3},D_{4}\},$ independent of $p$, such that 
$$[G:Z(G) \leq  |{\rm Aut}(S)| \leq (Dp)^{n-1}.$$ For later convenience, we note that we certainly have $D \geq 256$. 

\medskip
This deals with the case that $G/Z(G)$ is almost simple, and, in particular, we have now dealt with the case $M = G$. From now on,we suppose that $G/Z(G)$ is not almost simple.

\medskip
\noindent {\bf Remark:}  We note that if $p\in \{3,5 \}$,  we can do somewhat better in this analysis. For in that case the simple group $S$ is neither an alternating group nor a sporadic simple group, so that we could replace $D$ by $D^{\ast} = 24^{6}$ in these cases.

\medskip
Suppose next that $M = F(M)$. Notice that the minimality of $M$ ensures that $M/Z(M)$ is an elementary Abelian $q$-group $Q$ for some prime $q$, 
and that $G$ acts faithfully and irreducibly on $Q$. By the minimal choice of $M$, we know that $M = Z(G)U$,
where $U = $ is either extra-special, (or the central product of a cyclic group of order $4$ with an extraspecial $2$-group), say of  order $q^{2d+1}$, and $G/MC_{G}(M)$ is isomorphic to a subgroup of ${\rm Sp}(2d,q)$, where $q^{d}$ is the dimension of a faithful irreducible complex representation  of $M$.

\medskip
Then $$[G:MC_{G}(M)] \leq q^{2d^{2}+d}$$ and $[G:C_{G}(M)] \leq q^{2d^{2}+3d}.$ By Lemma 4.4 of [8], we
have $$q^{2d^{2}+ 3d} \leq \frac{16^{q^{d}}}{q^{d}}.$$ This is certainly at most $16^{q^{d}-1}$ when $q^{d} \geq 16$.
Hence we certainly have $$q^{2d^{2}+3d} \leq 256^{q^{d}-1}$$ for $q$ prime and $d$ any  positive integer.
Hence in this situation,we have $$[G:C_{G}(M)] \leq 256^{q^{d}-1},$$ where $q^{d} > 1$ is a prime power divisor of $n.$

\medskip
If we also have $C_{G}(M) \leq M$ when $M=F(M)$, then we certainly have $$[G:F(G)]  \leq D^{n-1}.$$ 

\medskip
We now prove by induction on the degree  $n$ that (whatever the structure of $M$), we have $[G:F(G)] \leq(Dp)^{n-1}$, where $D$ is the constant defined above.
this being true when $n =1$, and in also in the case that $F^{\ast}(G) = Z(G)L$, where $L$ is quasi-simple. For inductive purposes, we assume that the above constant $D$ works
 (for any choice of prime $p$) for all finite complex linear groups $H$ of degree less than $n$. That is, we assume that we have $$[H:F(H)]< (Dp)^{m-1}$$ for whenever $m < n$ 
and $p$ is a prime such that $H$ is a finite $p^{\prime}$-subgroup of ${\rm GL}(m,\mathbb{C}).$

\medskip
In the case where our minimally chosen non-central normal subgroup $M$ satisfies $C_{G}(M) = Z(M)$ and $M =F(M)$, we have demonstrated above that 
$[G:Z(G)]  \leq D^{n-1}$ since $D  \geq 256$ and $n \geq q^{d}.$

\medskip
Suppose next that $M = Z(G)E(M).$ 

\medskip
In this case, $E = E(M) = L_{1}L_{2} \ldots L_{t}$, where the $L_{i}$ are $G$-conjugate components of $G$, all with common centre $Z(L_{1}) \leq Z(E) \leq Z(G).$
Then $M$ is isomorphic to subgroup of ${\rm GL}(r^{t},\mathbb{C})$ where $r >1$ is a divisor of $n$, and each $L_{i}$ is isomorphic to an irreducible subgroup
of ${\rm GL}(r,\mathbb{C}).$ 

\medskip
We have already dealt with the case $t = 1$ and $C_{G}(M) = Z(G),$ so we exclude this case from further consideration. Suppose then that $C_{G}(M) = Z(G)$ and that $t >1$. 

\medskip
Notice that $L_{1}$ is now isomorphic to a subgroup of ${\rm GL}(r, \mathbb{C}).$ Let $L_{1}/Z(L_{1}) \cong S$, a non-Abelian simple group.
By the analysis used in the case that $G/Z(G)$ is almost simple, we know that  $|S| < (Dp)^{r-1}$ for the constant $D$ defined above.

\medskip
Now  $[G:Z(G)]  \leq |{\rm Aut}(S) \wr X|$, where $X$ is a (transitive) $p^{\prime}$-subgroup of $Sym_{t}$. Then we have $$[G:Z(G)] \leq (Dp)^{rt-t}p^{t-1} \leq (Dp)^{rt-1} \leq (Dp)^{r^{t}-1},$$ so certainly $$[G:Z(G)] \leq (Dp)^{n-1}.$$

\medskip
Hence in all cases in which $C_{G}(M) \leq M$, (including the case $M = G$), we have $$C_{G}(M)=Z(M) = Z(G) = F(G),$$ and we have demonstrated that $$[G:F(G)] \leq(Dp)^{n-1}.$$ 

\medskip
Now we suppose that $H = C_{G}(M) \not \leq M$, in which case $M \neq G.$ 

\medskip
In the case $M = Z(G)E(M)$, we see that $H$ is isomorphic to a subgroup of ${\rm GL}(\frac{n}{r^{t}}, \mathbb{C}).$ 
This time, (keeping the notation above, and by the inductive hypothesis), we have $$[G:H]  \leq |{\rm Aut}(S) \wr X|,$$ where $X$ is a (transitive) $p^{\prime}$-subgroup of the symmetric  group $Sym_{t}$. Then we have $$[G:H] \leq (Dp)^{rt-t}p^{t-1} \leq (Dp)^{rt-1} \leq (Dp)^{r^{t}-1}.$$ 

\medskip
Also, by the inductive hypothesis, we also have $$[H:F(H)] \leq (Dp)^{\frac{n}{r^{t}}-1}$$. Hence we have $$[G:F(H)] \leq (Dp)^{r^{t}-1+\frac{n}{r^{t}}-1} .$$
Now we have $ x + \frac{n}{x}\leq 2 +\frac{n}{2}$ whenever $x > 1$ is an integer which divides $n$,  and $r^{t} \geq 2$ is a divisor of $n$, so that $$[G:F(H)] \leq 
(Dp)^{\frac{n}{2}} \leq (Dp)^{n-1}$$ and certainly $$[G:F(G)] \leq(Dp)^{n-1}.$$

\medskip
Now we return to the case $M = F(M)$, but this time with $$H = C_{G}(M) \neq Z(M).$$  Now (using the earlier notation), we know that $H$ is isomorphic to a subgroup of 
${\rm GL}(\frac{n}{q^{d}},\mathbb{C})$, and by induction, we may suppose that  $$[H:F(H)] \leq (Dp)^{\frac{n}{q^{d}} -1}$$, 
while we still still have  $$[G:H] \leq 256^{q^{d}-1}.$$ 

\medskip
Hence we certainly have $$[G:F(H)] \leq (Dp)^{\frac{n}{q^{d}} -1} \times  D^{q^{d}-1}.$$ 
Now $x+ \frac{n}{x}$ decreases monotonically for $1 \leq x \leq \sqrt{n}$ and increases monotonically for $ x \geq \sqrt{n}.$  Since $q^{d}\geq 2$ and $q^{d}$ is a divisor of $n$, we conclude that 
$$[G:F(H)]  \leq D^{\frac{n}{2}} p^{\frac{n}{2} -1} < (Dp)^{n-1}.$$ Then we certainly have $[G:F(G)] \leq (Dp)^{n-1}$ in this case.

\medskip
Now we have proved that we always have $$[G:F(G)]  \leq (Dp)^{n-1}, $$ as required to complete the proof of Theorem A.

\medskip
We now complete the proof of Theorem B, given the existence of the constant $C$ in Theorem A. Let $p$ be a prime, $n$ be a positive integer, and let 
$G$ be a $p$-closed finite subgroup of ${\rm GL}(n,\mathbb{C}) $ with Sylow $p$-subgroup $P$. By the Schur-Zassenhaus Theorem (see, for example, Isaacs [7]),we have $G = HP$ for some  $p^{\prime}$-subgroup $H$ of $G$. Applying Theorem A to $H$, we see that there is an Abelian normal subgroup $A$ of $H$ with $[H:A] < (Cp)^{n-1}$.

\medskip
Notice that $K = PA$ is now a normal subgroup of $G$. Applying Dornhoff's theorem on solvable groups, there is a fixed constant $c$ (independent both of $n$ and of $p$) 
such that $[K:F(K)] < c^{n-1}.$  Now $F(K)$ is characteristic in the normal subgroup $K$ of $G$, so that $F(K) \lhd G.$ Hence certainly $[G:F(G)] \leq (cCp)^{n-1},$
and the existence of Abelian normal subgroup $A^{\ast}$ and the constant $C^{\ast}$ in Theorem B follows as in the proof of Theorem A.

\medskip
Suppose further that $P$ is Abelian (as well as normal in $G$). Then $PA^{\ast} \lhd G$, so that $PA^{\ast} \leq F(G)$ since $P$ and $A^{\ast}$ are both nilpotent.
However, $A^{\ast}$ has and Abelian Sylow $q$-subgroup for each $q \neq p$, while $P$ is an Abelian Sylow $p$-subgroup of $G$. Thus $PA^{\ast}$ is nilpotent with Abeilan Sylow 
subgroups for every prime $r$, and $PA^{\ast}$ is Abelian, while we certainly have $$[G:PA^{\ast}] \leq [G:A^{\ast}],$$so we may replace $A^{\ast}$ by $PA^{\ast}$ if necessary to ensure that $[G:A^{\ast}]$ is coprime to $p$.

\section{References}

\medskip
\noindent [1] Blichfeldt,H.,  \emph{Finite Collineation Groups},University of Chicago Press, Chicago, (1917).

\medskip
\noindent [2] Collins, M.J., \emph{On Jordan's Theorem for Complex Linear Groups}, Journal of Group Theory, 10 , (2007), 411-423.

\medskip
\noindent [3] Coulembier, K.,  Etingof, P. and Ostrik, V., \emph{Asymptotic Properties of Tensor Powers in Symmetric Tensor Categories}, to appear.

\medskip
\noindent [4] Conway, J.H., Curtis, R.T., Norton, S.P.,Parker, R.A and  Wilson, R.A., \emph{Atlas of Finite Groups}, Clarendon, Oxford, (1985).

\medskip
\noindent [5] Dornhoff, L., \emph{Group Representation Theory, Part A : Ordinary Representation Theory}, Marcel Dekker, New York, (1972).

\medskip
\noindent [6] Gow, R., \emph{On the number of characters in a $p$-block of a $p$-solvable group},  Journal of Algebra, 65, 2, (1980), 421-426.

\medskip
\noindent [7] Isaacs, I.M.,\emph{Character Theory of Finite Groups}, Academic Press, New York, (1976).

\medskip
\noindent [8] Kov\'acs,  L.G. and  Robinson, G.R., \emph{On the Number of Conjugacy Classes of a Finite Group}, Journal of Algebra, 160, (1993), 441-460.

\medskip
\noindent [9] P\'alfy, P. and Pyber, L.,\emph{Small Groups of Automorphisms}, Bulletin of the London Mathematical Society, 30, 4, (1998), 386-390.

\medskip
\noindent [10] Robinson, G.R., \emph{On Linear Groups}, Journal of Algebra, 131, (1990), 527-534.

\medskip
\noindent [11] Robinson, G.R., \emph {Bounding the Size of Permutation Groups and Complex Linear Groups of Odd Order}, Journal of Algebra, 335, (2011), 163-170.

\medskip
\noindent [12] Speiser, A., \emph{ Die Theorie der Gruppen vom endlicher Ordnung}, 3rd edition, (1937) ,Berlin.

\medskip
\noindent [13] Wagner, A., \emph{An Observation on the Degrees of Projective  Representations of the Symmetric and Alternating Group over an arbitrary Field},  Archiv der Mathematik, 29, (1977), 583-589.

\end{document}